\documentclass[12pt]{article}
\usepackage{amssymb}
\usepackage[french]{babel} 
\usepackage{amsmath,amsfonts}
\usepackage[applemac]{inputenc}

\def\thebibliograph#1#2{\section*{{\normalsize \bf #2}}\list
   {[\arabic{enumi}]}{\settowidth\labelwidth{[#1]}\leftmargin\labelwidth
     \advance\leftmargin\labelsep
     \usecounter{enumi}}
     \def\newblock{\hskip .11em plus .33em minus -.07em}
     \sloppy
     \sfcode`\.=1000\relax}

\newcommand{\dist}{{\mathcal D}'({\mathbb{R}}^n)}

\newcommand{\test}{{\mathcal D}({\mathbb{R}}^n)}

\newcommand{\tombstone}{\hspace{0cm}\hspace*{\fill} \rule{3mm}{3mm}
\\[2mm]}
\newcommand{\euclid}{{\mathbb{R}}^n}
\newcommand{\re}{{\mathbb{R}}}
\newcommand{\R}{{\mathbb{R}}^n}

\newcommand{\besovn}[3]{B_{{#2},{#3}}^{#1}({\mathbb R}^n)}

 \newcommand{\lizorn}[3]{F_{{#2},{#3}}^{#1}({\mathbb R}^n)}

   \newcommand{\spacen}[3]{E_{{#2},{#3}}^{#1}({\mathbb R}^n)}

\newtheorem{thm}{Th\'eor\`eme}
\newtheorem{prop}{Proposition}
\newtheorem{DE}{D\'efinition}

\newtheorem{LEM}{Lemme}

 \title{Localisation uniforme des espaces de Besov et de Li\-zorkin-Triebel}

 \author{Salah Eddine Allaoui et G\'erard Bourdaud}

\begin{document}
 \maketitle
 \date{}
 
\begin{abstract}

On \'etablit des caract\'erisations intrins\`eques des versions localis\'ees-uniformes des espaces de Besov ${\besovn{s}{p}{q}}$, avec  $p,q\in [1,+\infty]$, et de Lizorkin-Triebel ${\lizorn{s}{p}{q}}$ avec $q\in [1,+\infty]$ et $p\in [1,+\infty[$, quel que soit le nombre r\'eel $s>0$. \\

We give intrinsic characterisations for the uniformly localized versions of the  Besov spaces ${\besovn{s}{p}{q}}$, where $p,q\in [1,+\infty]$, and of the Lizorkin-Triebel spaces ${\lizorn{s}{p}{q}}$, where $q\in [1,+\infty]$ and $p\in [1,+\infty[$, whatever be the real number $s>0$.

\end{abstract}

{\it Mots-cl\'es:} {Espaces de Besov, Espaces de Lizorkin-Triebel, Localisation uniforme.}


\section{Introduction}

\`A tout espace norm\'e $E$ de fonctions sur $\euclid$, il est possible d'associer sa version localis\'ee-uniforme. Il s'agit de l'espace, not\'e
$E_{lu}$, des fonctions $f$ telles que
\begin{equation}\label{deflu} \sup_{a\in \euclid} \| (\tau_a\varphi)f\|_E< +\infty\,.\end{equation}
Ici $\tau_a$ d\'esigne l'op\'erateur de translation, d\'efini par $\tau_af(x):= f(x-a)$, et $\varphi$ est une fonction $C^\infty$ \`a support compact, positive, non nulle. Sous des hypoth\`eses standard, rappel\'ees au paragraphe \ref{local}, on peut montrer que $E_{lu}$ ne d\'epend aucunement du choix de la fonction auxiliaire $\varphi$.\\

Les espaces localis\'es-uniformes jouent un r\^ole dans diverses questions d'analyse math\'e\-ma\-tique. Par exemple, si $E$ est une alg\`ebre de Banach de fonctions, pour la multiplication usuelle, il est naturel de conjecturer que l'ensemble des multiplicateurs de $E$ est pr\'ecis\'ement $E_{lu}$ --- conjecture confirm\'ee dans le cas des espaces de Sobolev $H^s(\euclid)$ pour $s>n/2$, voir \cite[p.~58]{Bou_95} et \cite[p.~151]{P}. Ils interviennent aussi dans la caract\'erisation des fonctions qui op\`erent, par composition \`a gauche, sur certains espaces de fonctions. Ainsi, les fonctions de $\re$ dans $\re$ qui op\`erent, en ce sens, sur l'espace de Sobolev critique $W^m_p(\euclid)$, o\`u l'entier $m$ v\'erifie $m=n/p>1$, sont pr\'ecis\'ement celles dont les d\'eriv\'ees appartiennent localement-uniform\'ement \`a $W^{m-1}_p(\re)$ \cite{BouIn}. L'extension de ce th\'eor\`eme aux espaces de Sobolev fractionnaires est une question ouverte. Cependant on a pu \'etablir un r\'esultat partiel \cite{GB-SA}: si une fonction op\`ere sur l'espace de Besov $B^s_{p,q}(\euclid)$, avec $s=n/p>1$ et $q>1$, alors sa d\'eriv\'ee appartient localement-uniform\'ement \`a $B^{s-1}_{p,q}(\re)$.\\

Il semble d\`es lors pertinent de d\'ecrire les localisations uniformes des espaces de Sobolev fractionnaires de fa\c con intrins\`eque, c'est-\`a-dire sans recourir \`a une fonction auxiliaire telle que la fonction $\varphi$ utilis\'ee dans (\ref{deflu}). \`A cet \'egard, rappelons que, pour les espaces $L_p(\re^n)_{lu}$, une telle description est bien connue, voir la proposition \ref{Lplu}.
 C'est ce type de description, \`a l'aide d'int\'egrales portant sur les translat\'ees d'une boule fixe, que nous mettrons en \'evidence pour
les localisations uniformes des espaces de Sobolev fractionnaires.

\subsection*{Plan} La section \ref{local} sera d\'evolue \`a des g\'en\'eralit\'es sur la localisation uniforme.
Dans la section \ref{DBLT}, on rappellera les d\'efinitions des espaces Besov et de Lizorkin-Triebel  et on \'enoncera les th\'eor\`emes principaux, qui seront \'etablis dans la derni\`ere section.
\subsection*{Notations et rappels}

La norme d'une fonction $f$ dans l'espace de Lebesgue $L_p(\R)$ est not\'ee $\|f\|_p$.
Comme d'habitude $c, c_1,...$ d\'esignera une constante positive pouvant d\'ependre de $n, s, p, q$ et des fonctions auxiliaires; sa valeur pourra changer d'une occurrence \`a l'autre.
Rappelons une in\'egalit\'e classique (voir, par exemple, \cite[II.20, p.~44]{Bou_95}):

\begin{LEM}\label{MALPHA} Pour tout $q\in [1,+\infty[$ et tout r\'eel $\alpha$, il existe $c>0$ tel que
\[ \sup _{0<t\leq 1/2}  t^\alpha u(t)\leq c\left( \int_0^1 \left(t^\alpha u(t)\right)^q \frac{{\rm d}t}{t}\right)^{1/q}\]
pour toute fonction croissante $u$ sur l'intervalle $]0,1]$.
\end{LEM}

\section{G\'en\'eralit\'es sur la localisation uniforme}\label{local}

Un espace de Banach de distributions (E.B.D.) sur $\R$ est un sous-espace vectoriel $E$ de $\mathcal{D'}(\R)$ muni d'une norme compl\`ete $\|-\|_E$ telle que l'injection canonique $E\hookrightarrow\mathcal{D'}(\R)$ soit continue. On dit que l'espace $E$ est un $\mathcal{D}(\R)$-{\em module} si $\phi f\in E$ pour tout $\phi\in\mathcal{D}(\R)$ et tout $f\in E$.
Un E.B.D. $E$ est {\em isom\'etriquement invariant par translation} si $\tau_af\in E$ et $\|\tau_af\|_E=\|f\|_E$ pour tout $f\in E$ et tout $a\in\R$. 
\begin{prop}\label{luspace} Soit  $E$ un $\test$-module 
isom\'etriquement invariant par translation.  Pour toute distribution $f$, les deux propri\'et\'es suivantes sont \'equi\-va\-lentes:

(i) Il existe une fonction positive non nulle $\varphi\in \test$ v\'erifiant (\ref{deflu}).

(ii) Pour toute fonction $\varphi\in \test$, on a la propri\'et\'e (\ref{deflu}).
\end{prop}
{\bf Preuve.} Voir \cite[p.~57]{Bou_95}.

\tombstone

Si une distribution $f$ satisfait l'une des deux conditions \'equivalentes de la proposition \ref{luspace}, on dit que
 $f$ appartient
{\em localement uniform\'ement} \`a $E$; l'ensemble de ces distributions
est not\'e $E_{lu}$. Soit  une fonction positive non nulle $\varphi\in \test$.
On montre facilement que $E_{lu}$ est un E.B.D. pour la norme
\[ \|f\|_{E_{lu}}:=\sup_{a\in \R} \|(\tau_a\varphi)\,f\|_E \,.\]
 De la preuve de la proposition \ref{luspace}, il r\'esulte qu'\`a \'equivalence pr\`es, la norme de $E_{lu}$ ne d\'epend pas du choix de la fonction $\varphi$.\\

Si $E$ est un $E.B.D.$ et $m$ un entier positif, on peut
consid\'erer  {\em l'espace de Sobolev $W^m(E)$ d'ordre $m$} de base $E$,  \`a savoir
\[ W^m(E):= \{ f\in \dist\,:\, f^{(\alpha)} \in E\,\quad \mathrm{pour\,tout}\quad\, |\alpha|\leq m\}\,.\]
 $W^m(E)$ est un E.B.D. pour la norme
\[ \|f\|_{W^m(E)}:= \sum_{|\alpha|\leq m} \|f^{(\alpha)}\|_E\,.\]

\begin{prop}\label{Sobolev_m} Si $E$ est un $\test$-module, isom\'etriquement invariant par translation, il en est de m\^eme pour $W^m(E)$ et on a

\[ \left(W^m(E)\right)_{lu} = W^m\left(E_{lu}\right)\,,\]
avec des normes \'equivalentes.

\end{prop}

{\bf Preuve.}
Elle r\'esulte ais\'ement de la formule de Leibniz et de la proposition \ref{luspace}.

\tombstone

Nous terminerons cette section en rappelant la description de $L_p(\re^n)_{lu}$.
La preuve facile est laiss\'ee au lecteur.

\begin{prop}\label{Lplu}
Soit $p\in[1,+\infty[$.  Soit ${\mathbb B}$ une boule ouverte (ou un cube ouvert) dans $\re^n$. Alors une fonction mesurable $f$ sur $\R$ appartient \`a $L_p(\re^n)_{lu}$ si et seulement si

\[\sup_{a\in \re^n} \left( \int_{{\mathbb B}+a} |f(x)|^p{\rm d}x\right)^{1/p} <+\infty\,;\]
 de plus l'expression ci-dessus est \'equivalente \`a la norme $\left\|f\right\|_{L_p(\re^n)_{lu}}$.
 \end{prop}

%

\section{D\'efinitions des espaces fonctionnels et \'enonc\'es des th\'eor\`emes}\label{DBLT}

\`A toute fonction $f$, d\'efinie sur $\euclid$, et tout $h\in \euclid$, on associe la fonction
$\Delta_{h}f$, d\'efinie par $\Delta_{h}f:=\tau_{-h}f -f$.
Pour tout $p\in [1,+\infty]$, tout ensemble bor\'elien $A$ de $\R$, toute fonction mesurable $f$ sur $\R$ et tout $t>0$, on pose
\[ \omega_{p,A}(f,t):= \sup_{|h|\leq t} \left(\int_A \left|\Delta_{h}f(x)\right|^p {\rm d}x\right)^{1/p}\,,\]
\[ \eta_{p,A}(f,t):= \sup_{|h|\leq t} \left(\int_A \left|\Delta^2_{h}f(x)\right|^p {\rm d}x \right)^{1/p}\,.\]
On note simplement $\omega_{p}:=\omega_{p,\R}$, de m\^eme pour $\eta.$

\begin{DE} Soient  $0<s<1$, $p,q\in [1,+\infty]$. L'espace $\besovn{s}{p}{q}$ est l'ensemble des fonctions $f$ v\'erifiant
\[ \|f\|_{\besovn{s}{p}{q}}:= \|f\|_p+\left(\int_0^{1}(t^{-s}\omega_{p}(f,t))^q \frac{{\rm d}t}{t}\right)^{1/q}  < + \infty\,.\]
\end{DE}
\begin{DE} Soient  $0<s<1$, $q\in [1,+\infty]$, $1\leq p<\infty$.
L'espace de Lizorkin-Triebel
$\lizorn{s}{p}{q}$ est l'ensemble des fonctions $f$ v\'erifiant
\[ \|f\|_{\lizorn{s}{p}{q}}:= \|f\|_p+\left( \left(\int^1_{0}\left(t^{-s-n}\int_{|h|\leq t} |\Delta_{h}f(x)|^q\,{\rm d}h\right)^q\frac{{\rm d}t}{t}\right)^{p/q}\,{\rm d}x\right)^{1/p} < + \infty\,.\]
\end{DE}

Rappelons qu'on obtient les m\^emes espaces fonctionnels, avec des normes \'equivalentes, en rempla\c cant, dans les d\'efinitions pr\'ec\'edentes, l'int\'egrale $\int_0^1$ par l'int\'egrale $\int_0^r$, o\`u $r$ est n'importe quel r\'eel positif fix\'e.
Quand il n'y a pas lieu de distinguer entre les deux types d'espaces, $B$ ou $F$,
nous posons $\spacen{s}{p}{q} :=\besovn{s}{p}{q}$ ou $\lizorn{s}{p}{q}$.
Les espaces d'ordre $1$, c'est-\`a-dire les espaces $\spacen{1}{p}{q}$, se d\'efinissent de la m\^eme fa\c con,
{\em \`a condition de remplacer} $\Delta_h$ par l'op\'erateur de diff\'erence seconde  $\Delta_h^2$ (et donc, dans le cas Besov, $\omega$ 
par $\eta$). Les espaces d'ordre sup\'erieur \`a $1$ sont, par d\'efinition, les espaces de Sobolev bas\'es sur les espaces d'ordres compris entre $0$ et $1$:

\begin{DE}\label{ppr4}
Soient $s>1$ et  $m$ l'entier tel que $m<s\leq m+1.$ Alors $\spacen{s}{p}{q}$ est l'ensemble des  fonctions $f$ telles que  $f^{(\alpha)}\in\spacen{s-m}{p}{q}$ pour tout $\left|\alpha\right|\leq m$. Cet espace est muni de la norme
 \[\sum_{ |\alpha|\leq m} \|f^{(\alpha)}\|_{\spacen{s-m}{p}{q}}\,.\] 
  \end{DE}
 

Venons-en \`a la description intrins\`eque des espaces localis\'es-uniformes. On se limitera au cas $0<s\leq 1$, puisqu'il suffit
d'appliquer la proposition \ref{Sobolev_m} pour obtenir le cas g\'en\'eral.
Dans les \'enonc\'es suivants, ${\mathbb B}$ d\'esignera une boule (ou un cube) fix\'e de $\euclid$. On supposera $p,q\in [1,+\infty]$ ($p<\infty$ dans le cas Lizorkin-Triebel).
\begin{thm}\label{besovlu} Si $0<s<1$, alors $\besovn{s}{p}{q}_{lu}$ est l'ensemble des fonctions $f$
telles que
\begin{equation}\label{Besovlu}\sup_{a\in \R} \left(\int_0^1 \left( t^{-s} \omega_{p,{\mathbb B}+a}(f,t)\right)^q\frac{{\rm d}t}{t}\right)^{1/q}+\left\|f\right\|_{L_p(\R)_{lu}}<+\infty\,.\end{equation}
De plus l'expression ci-dessus 
est une norme \'equivalente sur $\besovn{s}{p}{q}_{lu}$.
\end{thm}

\begin{thm}\label{besovlus=1}  $\besovn{1}{p}{q}_{lu}$ est l'ensemble des fonctions $f$ telles que
\begin{equation}\label{Besovlus=1}\sup_{a\in \R} \left(\int_0^1 \left( t^{-1} \eta_{p,{\mathbb B}+a}(f,t)\right)^q\frac{{\rm d}t}{t}\right)^{1/q}+\left\|f\right\|_{L_p(\R)_{lu}}<+\infty\,.\end{equation}
De plus l'expression ci-dessus 
est une norme \'equivalente sur $\besovn{1}{p}{q}_{lu}$.
\end{thm}

\begin{thm}\label{Lizornlu} Si $0<s<1$,  alors $\lizorn{s}{p}{q}_{lu}$ est l'ensemble des fonctions $f$ telles que
\begin{equation}\label{LLizornlu}\sup_{a\in \R}\left\|\left(\int_0^1 \left( t^{-s-n}\int_{\left|h\right|\leq t}\left|\Delta_{h}f(.)\right|{\rm d}h\right)^q\frac{{\rm d}t}{t}\right)^{1/q}\right\|_{L_p({\mathbb B}+a)}+\left\|f\right\|_{L_p(\R)_{lu}}<+\infty\,.\end{equation}
De plus l'expression ci-dessus 
est une norme \'equivalente sur $\lizorn{s}{p}{q}_{lu}$.
\end{thm}

\begin{thm}\label{Lizornlus=1}  $\lizorn{1}{p}{q}_{lu}$ est l'ensemble des fonctions $f$ telles que
\begin{equation}\label{LLizornlus=1}\sup_{a\in \R}\left\|\left(\int_0^1 \left( t^{-n-1}\int_{\left|h\right|\leq t}\left|\Delta^2_{h}f(.)\right|{\rm d}h\right)^q\frac{{\rm d}t}{t}\right)^{1/q}\right\|_{L_p({\mathbb B}+a)}+\left\|f\right\|_{L_p(\R)_{lu}}<+\infty\,.\end{equation}
De plus l'expression ci-dessus
est une norme \'equivalente sur $\lizorn{1}{p}{q}_{lu}$.
\end{thm}

 \section{Preuves des th\'eor\`emes}
 Sans perte de g\'en\'eralit\'e, on peut supposer que ${\mathbb B}$ est la boule unit\'e de $\R$.
 Dans cette section,
 on fixe deux fonctions $\varphi_0$ et $\varphi_1$ dans $\test$, telles que:
 \begin{itemize}
 \item $0\leq \varphi_0\leq 1$, $\varphi_0$ est non nulle et port\'ee par 
 ${\mathbb B}/4$, 
 \item $\varphi_1(x)=1$ sur $4{\mathbb B}$.
 \end{itemize}

\subsection {Preuve du th\'eor\`eme \ref{besovlu}}

On utilisera la formule suivante, valable pour tout $h\in\R$ et  toutes fonctions $f$ et $g$ sur $\R$:
\begin{equation}\label{1diff}
\Delta_h(fg) = (\Delta_hf)(\tau_{-h}g)
 + f(\Delta_hg)\,.\end{equation}
D\'esignons par $A(f)$ le premier terme de l'in\'egalit\'e (\ref{Besovlu}).

\subsubsection{\'Etape 1}

Soit $f$ une fonction telle que $A(f)<\infty$.
Par la formule (\ref{1diff}), on a, pour tous $a,h\in \R$ et $|h|\leq t\leq 1/2$,
\[ \left(\int_{\R} \left| \Delta_{h}( (\tau_a\varphi_0) f) (x)\right|^p{\rm d}x\right)^{1/p}\]
  \[\leq \left(\int_{\R}  \left| \Delta_{h}f(x)\varphi_0(x+h-a)\right|^p{\rm d}x\right)^{1/p}+
    \left(\int_{\R} |f(x)|^p\left| \Delta_{h} (\tau_a\varphi_0) (x)\right|^p{\rm d}x\right)^{1/p}\]
    \[\leq  \left(\int_{{\mathbb B}+a}   \left| \Delta_{h}f (x)\right|^p{\rm d}x\right)^{1/p} +
    t\left\|\nabla\varphi_0\right\|_\infty\, \left(\int_{{\mathbb B}+a}   \left|f (x)\right|^p{\rm d}x\right)^{1/p}\]
      \[\leq c_1\,\left( \omega_{p,{\mathbb B}+a}(f,t) +
    \,t\,  \|f\|_{L^p(\R)_{lu}}\right)\,.\]
Par la condition $s<1$, on voit  que
\[  \left(\int_0^{1/2} \left( t^{-s} \omega_{p}((\tau_a\varphi_0)f,t)\right)^q\frac{{\rm d}t}{t}\right)^{1/q} \]
\[\leq c_1  \left(\left(\int_0^{1/2} \left( t^{-s} \omega_{p,{\mathbb B}+a}(f,t)\right)^q\frac{{\rm d}t}{t}\right)^{1/q}
+ \left(\int_0^{1/2} \left( t^{1-s} \right)^q\frac{{\rm d}t}{t}\right)^{1/q}  \|f\|_{L^p(\R)_{lu}}\right).\]
L'expression ci-dessus \'etant major\'ee par $c_2A(f)$, pour une certaine constante $c_2$,
il vient 
\[ \sup_{a\in \R} \|(\tau_a\varphi_0)f \|_{ \besovn{s}{p}{q}}\leq c_3 A(f)\,. \]

\subsubsection{\'Etape 2}

Soit $f\in\besovn{s}{p}{q}_{lu}$. On voit aussit\^ot que
$ \Delta_{h}((\tau_a\varphi_1)f)(x) = \Delta_{h}f(x)$ pour tout $a\in \R$,  tout $x\in {\mathbb B}+a$, et tout $|h|\leq 1$.
On en d\'eduit ais\'ement que
\[ A(f) \leq c_4 \sup_{a\in \R} \|(\tau_a\varphi_1)f \|_{ \lizorn{s}{p}{q}}\,.\]

\subsection {Preuve du th\'eor\`eme \ref{besovlus=1}}

On d\'esignera par $A(f)$ le premier terme de l'in\'egalit\'e (\ref{Besovlus=1}) et on posera
\[ M_{p,a}(f):= \sup_{0<t\leq 1/2} \,\frac{1}{t}\eta_{p,{\mathbb B}+a}(f,t)\,.\]

 \subsubsection{R\'esultats pr\'eliminaires  }
On dispose des formules suivantes, o\`u $k\in\mathbb{N^*}$, $h\in\R$ et o\`u $f$ et $g$ sont des fonctions quelconques:
\begin{equation}\label{2diff}
\Delta^2_h(fg) = (\Delta^2_hf)(\tau_{-2h}g) + (\Delta^2_hg)(\tau_{-h}f)
 +(\Delta_h f)(\Delta_{2h}g)\,,\end{equation}
 \begin{equation}\label{3diff}
 \Delta_{h} = 2^{-k} \Delta_{2^kh} - \sum_{l=0}^{k-1} 2^{-l-1}  \Delta^2_{2^lh}\, .
 \end{equation}
La premi\`ere est imm\'ediate, la seconde s'obtient facilement par r\'ecurrence sur $k$.

\begin{LEM}\label{Marchaud} Il existe $c>0$ tel que
\[ \omega_{p,{\mathbb B}+a}(f,t)\leq
ct\left\{ \left(\int_{2{\mathbb B}+a} |f(x)|^p\,{\rm d}x\right)^{1/p} + M_{p,a}(f)\,|\ln t|
\right\}\,,\]
pour tout $0<t\leq 1/2,$ tout $a\in \R$ et toute fonction localement int\'egrable $f.$ 
\end{LEM}
{\bf Preuve.} Le lemme est une variante de l'in\'egalit\'e classique de Marchaud.
On d\'efinit l'entier $k\geq 1$ par l'encadrement $2^{-k-1}<t\leq 2^{-k}.$ 
De la formule (\ref{3diff}), on d\'eduit, pour $|h|\leq t$,
\[ \left(\int_{{\mathbb B}+a}  \left| \Delta_{h}f(x)\right|^p{\rm d}x\right)^{1/p}\]
\[\leq 2^{-k}\left(\int_{{\mathbb B}+a}  \left| \Delta_{2^kh}f(x)\right|^p{\rm d}x\right)^{1/p}
+ \sum_{l=0}^{k-1} 2^{-l-1} \left(\int_{{\mathbb B}+a}  \left| \Delta^2_{2^lh}f(x)\right|^p{\rm d}x\right)^{1/p}\]
\[\leq 2^{-k+1}\,\left(\int_{2{\mathbb B}+a}  |f(x)|^p\,{\rm d}x\right)^{1/p} +\sum_{l=0}^{k-1} 2^{-l-1}(2^{l-k}M_{p,a}(f))\,,\]
\[\leq 4t\,\left(\int_{2{\mathbb B}+a}  |f(x)|^p\,{\rm d}x\right)^{1/p} +\frac{1}{\ln2}t\left|\ln t\right|M_{p,a}(f)\,,\]
ce qui conclut la preuve du lemme \ref{Marchaud}.

\subsubsection{\'Etape 1}

Soit $f$ une fonction telle que $A(f)<\infty$. Par la
formule (\ref{2diff}), il vient, pour $\left|h\right|\leq t\leq1/4$,

\[ \left(\int_{\R} \left| \Delta^2_{h}( (\tau_a\varphi_0) f) (x)\right|^p{\rm d}x\right)^{1/p}\leq \left(\int_{\R}  \left| \Delta^2_{h}f(x)\varphi_0(x+2h-a)\right|^p{\rm d}x\right)^{1/p}\] \[+
    \left(\int_{\R} |f(x+h)|^p\left| \Delta^2_{h} (\tau_a\varphi_0) (x)\right|^p{\rm d}x\right)^{1/p}\]
    \[+ \left(\int_{\R}  \left| \Delta_{h}f(x)(\varphi_0(x+2h-a)- \varphi_0(x-a))\right|^p{\rm d}x\right)^{1/p}
 \]\[\leq \left(\int_{{\mathbb B}+a}   \left| \Delta^2_{h}f (x)\right|^p{\rm d}x\right)^{1/p} +
    c_1t^2 \left(\int_{{\mathbb B}+a}   \left|f (x)\right|^p{\rm d}x\right)^{1/p}
   +c_2t\left(\int_{{\mathbb B}+a}   \left| \Delta_{h}f (x)\right|^p{\rm d}x\right)^{1/p}\]
    \[\leq\,c_3( \eta_{p,{\mathbb B}+a}(f,t) +\,t^2\,  \|f\|_{L^p(\R)_{lu}}+\,t\,\omega_{p,{\mathbb B}+a}(f,t))\,,\]
et donc
\[  \left(\int_0^{1/4} \left( t^{-1} \eta_{p}((\tau_a\varphi_0)f,t)\right)^q\frac{{\rm d}t}{t}\right)^{1/q} \]
\[\leq c_3  \left(\int_0^{1/4} \left( t^{-1} \eta_{p,{\mathbb B}+a}(f,t)\right)^q\frac{{\rm d}t}{t}\right)^{1/q} +   c_3\|f\|_{L^p(\R)_{lu}}\, \left( \int_0^{1/4} t^{q-1}{\rm d}t\right)^{1/q}\]
\[+c_3\left(\int_0^{1/4} \left(\omega_{p,{\mathbb B}+a}(f,t)\right)^q\frac{{\rm d}t}{t}\right)^{1/q}\,.\]
En cons\'equence
\[ \sup_{a\in \R} \|(\tau_a\varphi_0)f \|_{ \besovn{1}{p}{q}}\leq c_4\left(  A(f) + 
\sup_{a\in \R}\left(\int_0^{1/4} \left(\omega_{p,{\mathbb B}+a}(f,t)\right)^q\frac{{\rm d}t}{t}\right)^{1/q}\right)\,.\]
Gr\^ace au lemme
 \ref{Marchaud}, on a, pour tout $a\in \R$,
\[ \left(\int_0^{1/4} \left(\omega_{p,{\mathbb B}+a}(f,t)\right)^q\frac{{\rm d}t}{t}\right)^{1/q}\]\[
\leq c_5\|f\|_{L^p(\R)_{lu}}\, \left( \int_0^{1/4} t^{q-1}{\rm d}t\right)^{1/q}
+  c_6M_{p,a}(f)\, \left( \int_0^{1/4} t^{q-1}\, |\ln t|^q\,{\rm d}t\right)^{1/q}\,.\]
En appliquant le lemme \ref{MALPHA} \`a la fonction croissante $t\mapsto \eta_{p,{\mathbb B}+a}(f,t)$, on conclut que
  \[\sup_{a\in \R} \|(\tau_a\varphi_0)f \|_{ \besovn{1}{p}{q}} \leq  c_7A(f)\,.\]

\subsubsection{\'Etape 2}

Soit $f\in\besovn{1}{p}{q}_{lu}$. En proc\'edant comme dans l'\'etape 2 de la preuve du th\'eor\`eme \ref{besovlu}, il vient
\[ A(f) \leq c_8 \|f\|_{\besovn{1}{p}{q}_{lu}}\,.\]

 \subsection {Preuve du th\'eor\`eme \ref{Lizornlu}}
 On d\'esignera par $A(f)$ le premier terme de l'in\'egalit\'e (\ref{LLizornlu}).

 \subsubsection{\'Etape 1} 

Soit $f$ une fonction telle que $A(f)<\infty$.
Par la formule (\ref{1diff}), nous avons
\[ \int_{\left|h\right|\leq t} \left| \Delta_{h}( (\tau_a\varphi_0) f) (x)\right|{\rm d}h\]
\[\leq\int_{\left|h\right|\leq t} \left| \Delta_{h}f (x)\right| \,\varphi_0(x+h-a)\,{\rm d}h+\left|f(x)\right|\int_{\left|h\right|\leq t} \left| \Delta_{h} (\tau_a\varphi_0)(x)\right|{\rm d}h.\]
On obtient
\[\left(\int_{\R}\left(\int_0^{1/2} \left( t^{-s-n}\int_{\left|h\right|\leq t}\left|\Delta_{h}( (\tau_a\varphi_0) f) (x)\right|{\rm d}h\right)^q\frac{{\rm d}t}{t}\right)^{p/q}{\rm d}x\right)^{1/p} \,\]
\[\leq \left(\int_{{\mathbb B}+a}\left(\int_0^{1/2} \left( t^{-s-n}\int_{\left|h\right|\leq t}\left|\Delta_{h}f (x)\right|{\rm d}h\right)^q\frac{{\rm d}t}{t}\right)^{p/q}{\rm d}x\right)^{1/p}\,+
c_1\left(\int_{{\mathbb B}+a}   \left|f (x)\right|^p{\rm d}x\right)^{1/p},\]
ce qui nous donne
\[ \sup_{a\in \R} \|(\tau_a\varphi_0)f \|_{ \lizorn{s}{p}{q}} \leq c_2 A(f)\,.\]

 \subsubsection{\'Etape 2} 
Supposons que $f\in\lizorn{s}{p}{q}_{lu}$. En proc\'edant comme dans l'\'etape 2 de la preuve du th\'eor\`eme \ref{besovlu}, il vient
\[ A(f)\leq c_3\, \|f \|_{ \lizorn{s}{p}{q}_{lu}}\,.\]
 
 \subsection {Preuve du th\'eor\`eme \ref{Lizornlus=1}}
On d\'esignera par $A(f)$ le premier terme de l'in\'egalit\'e (\ref{LLizornlus=1}). 

 \subsubsection{\'Etape 1} 

Soit $f$ une fonction telle que $A(f)<\infty$. Soit
\[G(x):=\left(\int_0^{1} \left( t^{-n-1}\int_{\left|h\right|\leq t}\left|\Delta^2_{h}f(x)\right|{\rm d}h\right)^q\frac{{\rm d}t}{t}\right)^{1/q}\,.\]
Par la
formule (\ref{2diff}), il vient, pour tous $a,x\in \R$ et $t>0$,
\[ \int_{\left|h\right|\leq t} \left| \Delta^2_{h}( (\tau_a\varphi_0) f) (x)\right|{\rm d}h\]
\[\leq\int_{\left|h\right|\leq t}  \left| \Delta^2_hf(x)\right| \,\varphi_0(x+2h-a)\,{\rm d}h+\int_{\left|h\right|\leq t} |f(x+h)|\left| \Delta^2_h (\tau_a\varphi_0) (x)\right|{\rm d}h\]
 \[+ \int_{\left|h\right|\leq t}  \left| \Delta_hf(x)\right| \left| \Delta_{2h}\tau_a\varphi_0(x)\right|{\rm d}h.\]
On en d\'eduit, pour tout $a\in\R$,
\[\left(\int_{\R}\left(\int_0^{1/16} \left( t^{-n-1}\int_{\left|h\right|\leq t}\left|\Delta^2_{h}( (\tau_a\varphi_0) f) (x)\right|{\rm d}h\right)^q\frac{{\rm d}t}{t}\right)^{p/q}{\rm d}x\right)^{1/p} \]
\[\leq\left(\int_{\frac {\mathbb B}2+a}\left(\int_0^{1/16} \left( t^{-n-1}\int_{\left|h\right|\leq t}\left|\Delta^2_{h}( (\tau_a\varphi_0) f) (x)\right|{\rm d}h\right)^q\frac{{\rm d}t}{t}\right)^{p/q}{\rm d}x\right)^{1/p}\]\[ \leq c_1(U(a)+V(a)+W(a))\,,\]
o\`u
\[U(a):=\left(\int_{\frac {\mathbb B}2+a}\left(\int_0^{1/16} \left( t^{-n-1}\int_{\left|h\right|\leq t}\left|\Delta^2_{h}f(x)\right|{\rm d}h\right)^q\frac{{\rm d}t}{t}\right)^{p/q}{\rm d}x\right)^{1/p} \,,\]
\[V(a):=\left(\int_{\frac {\mathbb B}2+a}\left(\int_0^{1/16} \left( t^{-n+1}\int_{\left|h\right|\leq t}\left|f(x+h)\right|{\rm d}h\right)^q\frac{{\rm d}t}{t}\right)^{p/q}{\rm d}x\right)^{1/p} \,,\]
\[W(a):=\left(\int_{\frac {\mathbb B}2+a}\left(\int_0^{1/16} \left( t^{-n}\int_{\left|h\right|\leq t}\left|\Delta_h f(x)\right|{\rm d}h\right)^q\frac{{\rm d}t}{t}\right)^{p/q}{\rm d}x\right)^{1/p} \,.\]
On voit facilement que
\begin{equation}\label{U+V} U(a)+V(a)\leq c_2\left(\int_{{\mathbb B}+a}G(x)^p{\rm d}x\right)^{1/p} +c_3\left\|f\right\|_{L_p({\mathbb B}+a) }
\,.\end{equation}
{\em Estimation de W(a).} Posons
$$
G_1(x):=\left(\int_0^{1/16} \left( t^{-n}\int_{\left|h\right|\leq t}\left|\Delta_h f(x)\right|{\rm d}h\right)^q\frac{{\rm d}t}{t}\right)^{1/q}.
$$
En d\'ecomposant l'intervalle $]0,1/16]$ en intervalles dyadiques et en utilisant des majorations \'evi-dentes, on obtient
$
G_1(x)\leq c_4G_2(x)$,
o\`u
$$
G_2(x):=\left(\sum_{j\geq4}\left(2^{jn}\int_{\left|h\right|\leq2^{-j}}\left|\Delta_h f(x)\right|{\rm d}h\right)^q\right)^{1/q}.
$$
Par le changement de variable $h'=2^{j-3}h,$ il vient
$$
G_2(x)=\left(\sum_{j\geq4}\left(\int_{\left|h\right|\leq1/8}\left|\Delta_{2^{-j+3}h} f(x)\right|{\rm d}h\right)^q\right)^{1/q}.
$$
Par (\ref{3diff}), on a
$$
\Delta_{2^{-j+3}h} = 2^{-j+3} \Delta_{h} - \sum_{\ell=0}^{j-4} 2^{-l-1}  \Delta^2_{2^{\ell-j+3}h}\,,
$$
d'o\`u
$
G_2(x)\leq c_5\left(G_3(x)+G_4(x)\right),
$
avec
$$
G_3(x):=\left(\sum_{j\geq4}\left(2^{-j}\int_{\left|h\right|\leq1/8}\left|\Delta_h f(x)\right|{\rm d}h\right)^q\right)^{1/q},
$$
et
\[
G_4(x):=\left(\sum_{j\geq4}\left(\int_{\left|h\right|\leq1/8}\sum^{j-4}_{\ell=0}2^{-\ell-1}\left|\Delta^2_{2^{\ell-j+3}h} f(x)\right|{\rm d}h\right)^q\right)^{1/q}.
\]
{\em Estimation de $G_3$.} On a aussit\^ot
$$
G_3(x)= c_6\int_{\left|h\right|\leq1/8}\left|\Delta_h f(x)\right|{\rm d}h.
$$
L'in\'egalit\'e de Minkowski nous donne, pour tout $a\in \R$,
$$
\left(\int_{\frac{{\mathbb B}}{2}+a}G_3(x)^p{\rm d}x\right)^{1/p}\leq\,c_6\int_{\left|h\right|\leq1/8}\left\{\int_{\frac{{\mathbb B}}{2}+a}\left|\Delta_h f(x)\right|^p{\rm d}x\right\}^{1/p}{\rm d}h
$$
$$
\leq\,c_6\int_{\left|h\right|\leq1/8}\left\{\int_{\frac{{\mathbb B}}{2}+a}\left|f(x+h)\right|^p{\rm d}x\right\}^{1/p}{\rm d}h+c_7\left(\int_{\frac{{\mathbb B}}{2}+a}\left|f(x)\right|^p{\rm d}x\right)^{1/p}
$$

$$
\leq\,c_6\int_{\left|h\right|\leq1/8}\left\{\int_{{\mathbb B}+a}\left|f(x)\right|^p{\rm d}x\right\}^{1/p}{\rm d}h+c_7\left(\int_{\frac{{\mathbb B}}{2}+a}\left|f(x)\right|^p{\rm d}x\right)^{1/p}
$$
$$
\leq c_8\left(\int_{{\mathbb B}+a}\left|f(x)\right|^p{\rm d}x\right)^{1/p}.
$$
{\em Estimation de $G_4$.}  Par le lemme \ref{MALPHA}, on a, pour tout $x\in \R$ et $0<t\leq 1/2$,
\begin{equation}\label{G4} 
\int_{\left|h\right|\leq t}\left|\Delta^2_{h} f(x)\right|{\rm d}h\leq c_9t^{n+1}G(x)\,.
\end{equation}
En raison du plongement $\ell_1\hookrightarrow \ell_q$, on a
$$
G_4(x)\leq \sum_{j\geq4}\int_{\left|h\right|\leq1/8}\sum^{j-4}_{\ell=0}2^{-\ell-1}\left|\Delta^2_{2^{\ell-j+3}h} f(x)\right|{\rm d}h.
$$
On v\'erifie facilement que
\[\int_{\left|h\right|\leq1/8}\left|\Delta^2_{2^{\ell-j+3}h} f(x)\right|{\rm d}h=2^{-3n}\,2^{(j-\ell)n}  \int_{\left|h\right|\leq2^{\ell-j}}\left|\Delta^2_{h} f(x)\right|{\rm d}h\,.\]
En combinant cette relation avec l'in\'egalit\'e (\ref{G4}),
on obtient $$
G_4(x)\leq c_{10}G(x)\sum_{j\geq4}\sum^{j-4}_{\ell=0}2^{-\ell-1}2^{\ell-j}=c_{11}G(x).
$$
Il vient donc, pour tout $a\in \R$,
$$
\left(\int_{\frac{{\mathbb B}}{2}+a}G_4(x)^p{\rm d}x\right)^{1/p}\leq c_{11}\left(\int_{{\mathbb B}+a}G(x)^p{\rm d}x\right)^{1/p}.
$$
En tenant compte des estimations obtenues pour $G_3$ et $G_4,$ on peut conclure que l'expression $W(a)$ est estim\'ee par
$$
\left(\int_{{\mathbb B}+a}G(x)^p{\rm d}x\right)^{1/p}+\left(\int_{{\mathbb B}+a}\left|f(x)\right|^p{\rm d}x\right)^{1/p}.
$$
En combinant avec (\ref{U+V}), on conclut que\[ \sup_{a\in \R} \|(\tau_a\varphi_0)f \|_{\lizorn{1}{p}{q}} \leq c_{12} A(f)\,.\]

 \subsubsection{\'Etape 2} 
Supposons que $f\in\lizorn{1}{p}{q}_{lu}$. En proc\'edant comme dans l'\'etape 2 de la preuve du th\'eor\`eme \ref{besovlu}, il vient
\[ A(f)\leq c_{13} \|f\|_{ F^1_{p,q}(\R)_{lu}}\,.\]

\newpage

Salah Eddine Allaoui\\
D\'epartement de Math\'ematique
 et Informatique\\
Universit\'e de Laghouat\\
Laghouat 03000\\
Alg\'erie\\
shallaoui@yahoo.fr\\

\vskip 2mm

G\'erard Bourdaud
 
 Universit\'e Paris Diderot, I.M.J. - P.R.G (UMR 7586)

B\^atiment Sophie Germain
       
       Case 7012

75205 Paris Cedex 13     \\                          
bourdaud@math.univ-paris-diderot.fr

\end{document}